\documentclass[]{article}

\usepackage{arxiv}

\usepackage[utf8]{inputenc} % allow utf-8 input
\usepackage[T1]{fontenc}    % use 8-bit T1 fonts
\usepackage{hyperref}       % hyperlinks
\usepackage{url}            % simple URL typesetting
\usepackage{booktabs}       % professional-quality tables
\usepackage{amsfonts}       % blackboard math symbols
\usepackage{nicefrac}       % compact symbols for 1/2, etc.
\usepackage{amsthm}
\usepackage{amsmath}
\theoremstyle{plain}
\newtheorem{theorem}{Theorem}[section]
\newtheorem{definition}{Definition}[section]
\newtheorem{proposition}{Proposition}[section]

\usepackage{blkarray}
\usepackage{multicol}

\title{Structure-preserving model order reduction of Hamiltonian systems}

\author{
 Jan S. Hesthaven \\
  Institute of Mathematics\\
  École Polytechnique Fédérale de Lausanne (EPFL)\\
  Lausanne, Switzerland \\
  \texttt{jan.hesthaven@epfl.ch} \\
   \And
 Cecilia Pagliantini \\
  Department of Mathematics and Computer Science\\
  Eindhoven University of Technology,\\
  Eindhoven, Netherlands \\
  \texttt{c.pagliantini@tue.nl} \\
  \And
 Nicolò Ripamonti \\
  Institute of Mathematics\\
  École Polytechnique Fédérale de Lausanne (EPFL)\\
  Lausanne, Switzerland \\
  \texttt{nicolo.ripamonti@epfl.ch} \\
}

\begin{document}
\maketitle
\begin{abstract}
We discuss the recent developments of projection-based model order reduction (MOR) techniques targeting Hamiltonian problems. Hamilton's principle completely characterizes many high-dimensional models in mathematical physics, resulting in rich geometric structures, with examples in fluid dynamics, quantum mechanics, optical systems, and epidemiological models. MOR reduces the computational burden associated with the approximation of complex systems by introducing low-dimensional surrogate models, enabling efficient multi-query numerical simulations.
However, standard reduction approaches do not guarantee the conservation of the delicate dynamics of Hamiltonian problems, resulting in reduced models plagued by instability or accuracy loss over time.
By approaching the reduction process from the geometric perspective of symplectic manifolds, the resulting reduced models inherit stability and conservation properties of the high-dimensional formulations.
We first introduce the general principles of symplectic geometry, including symplectic vector spaces, Darboux' theorem, and Hamiltonian vector fields. These notions are then used as a starting point to develop different structure-preserving reduced basis (RB) algorithms, including SVD-based approaches and greedy techniques. 
We conclude the review by addressing the reduction of problems that are not linearly reducible or in a non-canonical Hamiltonian form.  
\end{abstract}

% keywords can be removed
%\keywords{First keyword \and Second keyword \and More}

\section{Introduction}
The discretization of partial differential equations (PDEs) by classical methods like finite element, spectral method, or finite volume leads to dynamical models with very large state-space dimensions, typically of the order of millions of degrees of freedom to obtain an accurate solution. MOR \cite{schilders} is an effective method for reducing the complexity of such models while capturing the essential features of the system state. 
Starting from the Truncated Balanced Realization, introduced by Moore \cite{moore} in 1981, several other reduction techniques have been developed and flourished during the last 40 years, including the Hankel-norm 
reduction \cite{glover}, the proper orthogonal decomposition (POD) \cite{method_of_snapshots} and the Padé-via-Lanczos (PVL) algorithm \cite{feldmann}.
More recently, there has been a focus on the physical interpretability of the reduced models. Failure to preserve structures, invariants, and intrinsic properties of the approximate model, besides 
raising questions about the validity of the reduced models, has been associated with instabilities and exponential error growth, independently of the theoretical accuracy of the reduced solution space.
Stable reduced models have been recovered by enforcing constraints on the reduced dynamics obtained using standard reduction tools. Equality and inequality constraints have been considered to control the amplitude of the POD modes \cite{control_POD_modes}, the fluid temperature in a combustor \cite{huang}, and the aerodynamic coefficients \cite{zimmer}. 
Other methods directly incorporate the quantity of interest into the reduced system, producing \emph{inf-sup} stable \cite{infsup}, flux-preserving \cite{flux}, and skew-symmetric \cite{ripa} conservative reduced dynamics.
Even though great effort has been spent developing time integrators that preserve the symplectic flow underlying Hamiltonian systems, interest in geometric model order reduction initiated more recently, with efforts to preserve the Lagrangian structures \cite{lall_lagrangian}.\\
The remainder of the paper is organized as follows. In Section \ref{sec:general_symplectic_geometry}, we present the structure characterizing the dynamics of Hamiltonian systems and 
the concept of symplectic transformations. In Section \ref{sec:symplectic_reduced_methods}, we show that linear symplectic maps can be used to guarantee that the reduced models inherit the geometric formulation from the full dynamics. Different strategies to generate such maps are investigated in Section \ref{sec:PSD}, with thoughts on optimality results and computational complexities. A novel approach deviating from the linearity of the projection map is briefly discussed in Section \ref{sec:dynamical_low_rank}.
Finally, we discuss applications of structure-preserving reduction techniques to two more general classes of problems in Section \ref{sec:extension}, and some concluding remarks are offered in Section \ref{sec:conclusion}.

\section{Symplectic geometry and Hamiltonian systems}\label{sec:general_symplectic_geometry}
Let us first establish some definitions and properties concerning symplectic vector spaces.
\begin{definition}
Let $\mathcal{M}$ be a finite-dimensional real vector space and $\Omega:\mathcal{M}\times \mathcal{M}\mapsto \mathbb{R}$ a bilinear map. $\Omega$ is called anti-symmetric if
\begin{equation*}
    \Omega(u,v) = - \Omega(v,u), \qquad \forall u,v\in \mathcal{M}.
\end{equation*}
It is non-degenerate if
\begin{equation*}
    \Omega(u,v)=0, \forall u\in \mathcal{M}, \Rightarrow v=0.
\end{equation*}
\end{definition}
\begin{definition}
Let $\mathcal{M}$ be a finite-dimensional  vector space with $\Omega$ an anti-symmetric bilinear form on $\mathcal{M}$. The pair $(\mathcal{M},\Omega)$ is a symplectic linear vector space if $\Omega$ is non-degenerate. Moreover, $\mathcal{M}$ has to be $2n$-dimensional.
\end{definition}
Since we are interested in structure-preserving transformations, preserving the structure means to preserve the anti-symmetric bilinear form, as stated in the following definition.
\begin{definition}\label{def:symplectic_transform}
Let $(\mathcal{M}_1,\Omega_1)$ and $(\mathcal{M}_2,\Omega_2)$ be two symplectic vector spaces with $\text{dim}(\mathcal{M}_1)\geq\text{dim}(\mathcal{M}_2)$. The differentiable map $\varphi:\mathcal{M}_1\mapsto \mathcal{M}_2$ is called a symplectic transformation (symplectomorphism) if
\begin{equation*}
\varphi^{*}\Omega_2 = \Omega_1,    
\end{equation*}
where $\varphi^{*}\Omega_2$ is the pull-back of $\Omega_2$ 
with $\varphi$.
\end{definition}
One of the essential properties of Euclidean spaces is that all the Euclidean spaces of equal dimensions are isomorphic. For the symplectic vector spaces, a similar result holds, since two $2n$-dimensional symplectic vector spaces are symplectomorphic to one another. They therefore are fully characterized by their dimensions (As a consequence of the following theorem).  
\begin{theorem}[Linear Darboux' theorem \cite{darboux}]\label{thm:linear_darboux_theorem}
For any symplectic vector space $(\mathcal{M},\Omega)$, there exists a basis $\{ e_i, f_i \}_{i=1}^{n}$ of $\mathcal{M}$  such that
\begin{equation}\label{eq:properties_fomr_vector_space}
    \Omega(e_i,e_j)=0=\Omega(f_i,f_j), \qquad \Omega(e_i,f_j)=\delta_{ij}, \qquad \forall i,j=1,\dots,n.
\end{equation}
The basis is called Darboux' chart or canonical basis.
\end{theorem}

The proof of Theorem \ref{thm:linear_darboux_theorem} is based on a procedure similar to the Gram-Schmidt process to generate the symplectic basis, known as symplectic Gram-Schmidt \cite{symplectic_gram_schmidt}. \\
The canonical basis allows representing the symplectic form as
\begin{equation}\label{eq:fancy_symplectic_form}
\Omega(u,v)=\zeta^{\top}\mathbb{J}_{2n}\eta,
\end{equation}
where $\zeta,\eta\in\mathbb{R}^{2n}$ are the expansion coefficients of $u,v\in\mathcal{M}$ with respect to the basis $\{ e_i,f_i \}_{i=1}^{n}$ and
\begin{equation}\label{eq:Poisson_tensor}
    \mathbb{J}_{2n} = 
    \begin{bmatrix}
    \enspace 0_n & \enspace \mathbb{I}_n \enspace \\
    \enspace  -\mathbb{I}_n & \enspace 0_n \enspace
    \end{bmatrix}\in\mathbb{R}^{2n\times 2n},
\end{equation}
is known as the Poisson tensor, with $0_n\in\mathbb{R}^{n\times n}$ and $\mathbb{I}_n\in\mathbb{R}^{n\times n}$ denoting the zero and identity matrices, respectively. As a direct result, the matrix representation of the symplectic form $\Omega$ in the canonical basis is $\mathbb{J}_{2n}$. 
More generally, using a non-canonical basis, the form reduces to $\Omega(u,v)=\zeta^{\top}J_{2n}\eta$, with $J_{2n}$ being an invertible constant skew-symmetric matrix.\\
While symplectic vector spaces are helpful for the analysis of dynamical problems in Euclidean spaces and to define geometric reduced-order models, the constraint to the Euclidean setting is not generally adequate. In particular, the abstraction of the phase spaces of classical mechanics over arbitrary manifolds requires the definition of more general symplectic manifolds. We refer the reader to \cite{marsden} for a more comprehensive description of the topic. 
In this work, we limit ourselves to introducing a significant result regarding the evolution of the state of Hamiltonian systems.
\begin{definition}
Let $(\mathcal{M},\Omega)$ be a symplectic manifold and $H:\mathcal{M}\mapsto \mathbb{R}$ a $1$-form.
We refer to the unique vector field $\mathcal{X}_{H}$, which satisfies
\begin{equation*}
    i(\mathcal{X}_H)\Omega = \mathbf{d} H,
\end{equation*}
as the \emph{Hamiltonian vector field} related to $H$, where $i(\mathcal{X}_H)$ denotes the contraction operator and $\mathbf{d}$ is the exterior derivative. The function $H$ is called the \emph{Hamiltonian} of the vector field $\mathcal{X}_{H}$.
\end{definition}
Suppose $\mathcal{M}$ is also compact, then $\mathcal{X}_H$ is complete \cite{completeness_H_vector} and can be integrated, i.e., there exists an integral curve of $\mathcal{X}_H$, parametrized by the real variable $t$, that is the solution of
\begin{equation}\label{eq:Hamilton_equation}
    \dot{y}(t) = \mathcal{X}_H(y(t)).
\end{equation}
Equation \eqref{eq:Hamilton_equation} is referred to as Hamilton's equation of evolution or Hamiltonian system. Darboux' theorem, as a generalization of Theorem \ref{thm:linear_darboux_theorem}, states that two symplectic manifolds are only locally symplectomorphic.
Using this result, the Hamiltonian vector field $\mathcal{X}_H$ admits the local representation 
\begin{equation}\label{eq:local_vector_field}
    \mathcal{X}_H = \sum_{i=1}^{n} \dfrac{\partial H}{\partial f_i}\dfrac{\partial}{\partial e_i} - \dfrac{\partial H}{\partial e_i} \dfrac{\partial}{\partial f_i},
\end{equation}
with $\{ e_i, f_i \}_{i=1}^{n}$ is a local basis, leading to the following representation of \eqref{eq:Hamilton_equation}, expressed directly in terms of $H$.
\begin{proposition}
Let $(\mathcal{M},\Omega)$ be a $2n$-dimensional symplectic vector space and let $\{q_i,p_i\}_{i=1}^{n}$ be a canonical system of coordinates. Hamilton's equation is defined by
\begin{equation}\label{eq:Hamilton_equation_canonical}
\begin{cases}
\dfrac{dq_i}{dt} = \dfrac{\partial H}{\partial p_i},\\
\dfrac{dp_i}{dt} = -\dfrac{\partial H}{\partial q_i},
\end{cases}
\end{equation}
for $i=1,\dots, n$, which is a first order system in the $(q_i,p_i)$-space, or generalized phase-space.
\end{proposition}
Thus, if the state vector $y=(q_1,\dots,q_n,p_1,\dots,p_n)$ is introduced, \eqref{eq:Hamilton_equation_canonical} takes the form 
\begin{equation}\label{eq:Hamilton_equation_canonical_nabla}
    \dot{y}(t)=\mathbb{J}_{2n}\nabla_y H(y(t)),
\end{equation}
where $\nabla_y H$ is the naive gradient of $H$.
The flow of Hamilton's equation has some interesting properties.
\begin{proposition}\label{prop:flow_symplectic_transform}
Let $\phi_t$ be the flow of a Hamiltonian vector field $\mathcal{X}_H$. Then $\phi_t:\mathcal{M}\mapsto\mathcal{M}$ is a symplectic transformation.
\end{proposition}
We rely on a geometric perspective of linear vector spaces to highlight the importance of Proposition \ref{prop:flow_symplectic_transform}. 
Given two coefficient vectors $u$ and $v$ in $\mathbb{R}^{2n}$, the symplectic form \eqref{eq:fancy_symplectic_form} can be interpreted as the sum of the oriented areas of the orthogonal projection of the parallelogram defined by the two vectors on the $(q_i,p_i)$ planes.
Definition \ref{def:symplectic_transform}, in case of $2n$-dimensional symplectic vector space $(\mathcal{M},\Omega)$ with canonical coordinates, is equivalent to stating that a map $\phi:\mathbb{R}^{2n}\mapsto \mathbb{R}^{2n}$ is a symplectic transformation if and only if its Jacobian $\phi^{'}$ satisfies everywhere
\begin{equation}\label{eq:symplectic_condition_linear}
    (\phi^{'})^{\top}J_{2n}\phi^{'} = J_{2n}.
\end{equation}
Property \eqref{eq:symplectic_condition_linear} can be used to show that a symplectic transformation preserves the bilinear form $\Omega$ in the sense that \cite{marsden}
\begin{equation}
    \Omega(\phi(u),\phi(v)) = \Omega(u,v).
\end{equation}
Hence, the symplectic map $\phi$ represents a volume-preserving transformation. However, being symplectic is a more restrictive condition than being volume-preserving, as shown in the Non-squeezing theorem \cite{gromov}.\\
We conclude this section by noting that if the Hamiltonian function does not depend explicitly on time, its value is conserved along the solution trajectory.
\begin{proposition}\label{prop:Hamiltonian_conservation}
For Hamiltonian systems \eqref{eq:Hamilton_equation_canonical_nabla}, the Hamiltonian function is a first integral.
\end{proposition}

\section{Symplectic Galerkin projection}\label{sec:symplectic_reduced_methods}
The motivation of MOR is to reduce the computational complexity of dynamical systems in numerical simulations. In the context of structure-preserving projection-based reduction, two key ingredients are required to define a reduced model. First, we need a low-dimensional symplectic vector space that accurately represents the solution manifold of the original problem.  Then, we have to define a projection operator to map the symplectic flow of the Hamiltonian system onto the reduced space, while preserving its delicate properties. \\
Let us assume there exists a canonical basis $\{ e_i, f_i \}_{i=1}^{n}$ such that Hamilton's equation can be written in canonical form
\begin{equation}\label{eq:hamilton_equation_tbr}
    \begin{cases}
    \dot{y}(t)=\mathbb{J}_{2n}\nabla_y H(y(t)),  \\\
    y(0) = y_{0},
    \end{cases}
\end{equation}
and the related symplectic vector space is denoted by $(\mathcal{M},\Omega)$. Symplectic projection-based model order reduction adheres to the key idea of more general projection-based techniques \cite{hesthaven_rb} to approximate $y$ in a low-dimensional symplectic subspace $(\mathcal{A},\Omega)$ of dimension $2k$. In particular, we aim at $k\ll n$ to have a clear reduction, and therefore, significant gains in terms of computational efficiency. Let $\{\Tilde{e}_i,\Tilde{f}_i\}_{i=1}^{k}$ be a reduced basis for the approximate symplectic subspace and construct the linear map $\phi:\mathcal{A}\mapsto \mathcal{M}$ given by
\begin{equation}\label{eq:reduced_basis_approximation}
    y \approx \phi(z) = Az, 
\end{equation}
where
\begin{equation*}
    A = [\Tilde{e}_1,\dots , \Tilde{e}_k,\Tilde{f}_1,\dots, \Tilde{f}_k]\in\mathbb{R}^{2n\times 2k}.
\end{equation*}
$A$ belongs to the set of symplectic matrices of dimension $2n\times 2k$, also known as the symplectic Stiefel manifold, defined by
\begin{equation*}
    Sp(2k,\mathbb{R}^{2n}) := \{ L\in\mathbb{R}^{2n\times 2k}: L^{\top}J_{2n}L = J_{2k} \}.
\end{equation*}
Differential maps are often used to transfer structures from well-defined spaces to unknown manifolds. In this context, using the symplecticity of $A$, it is possible to show \cite{babak_symplectic} that Definition \ref{def:symplectic_transform} holds, with the right inverse of $\phi$ represented by $A$, and that there exists a symplectic form on $\mathcal{A}$ given by
\begin{equation}\label{eq:symplectic_form_subspace}
    \Tilde{\Omega} = \phi^{*}\Omega = A^{\top} \mathbb{J}_{2n} A = \mathbb{J}_{2k}. 
\end{equation}
As a result, $(\mathcal{A},\Tilde{\Omega})$ is a symplectic  vector space. In the following, for the sake of notation, we use $\mathcal{A}$ to indicate the reduced symplectic manifold paired with its bilinear form. \\
Given a symplectic matrix $A\in\mathbb{R}^{2n\times 2k}$, its symplectic inverse is defined as 
\begin{equation}
    A^{+} = \mathbb{J}_{2k}^{\top} A^{\top} \mathbb{J}_{2n}.
\end{equation}
Even though different from the pseudo-inverse matrix $(A^{\top}A)^{-1}A^{\top}$, the symplectic inverse $A^{+}$ plays a similar role and, in the following proposition, we outline its main properties \cite{peng_symplectic}.
\begin{proposition}\label{prop:properties_symplectic_inverse}
Suppose $A\in\mathbb{R}^{2n\times 2k}$ is a symplectic matrix and $A^{+}$ is its symplectic inverse. Then
\begin{multicols}{2}
\begin{itemize}
    \item $A^{+}A =\mathbb{I}_{2n}$.
    \item $(((A^{+})^{\top})^{+})^{\top}=A$.
    \item $(A^{+})^{\top}\in Sp(2k,\mathbb{R}^{2n})$.
    \item If A is  orthogonal then $A^{+}=A^{\top}$.
\end{itemize}
\end{multicols}
\end{proposition}
\noindent Using \eqref{eq:symplectic_form_subspace}, the definition of $A^{+}$ and the symplectic Gram-Schmidt process, it is possible to construct a projection operator 
    $\mathcal{P}_{\mathcal{A}} = A \mathbb{J}_{2k}^{\top} A^{\top} \mathbb{J}_{2n} = A A^{+}$,
that, differently from the POD orthogonal projection \cite{galerkin_orthogonal_projection}, can be used to approximate \eqref{eq:hamilton_equation_tbr} with  Hamiltonian system of reduced-dimension $2k$, characterized by the Hamiltonian function 
\begin{equation}\label{eq:reduced_Hamiltonian}
    H_{RB}(z) = H(Az).
\end{equation}
In particular, in the framework of Galerkin projection, using \eqref{eq:reduced_basis_approximation} in \eqref{eq:hamilton_equation_tbr}
yields
\begin{equation}\label{eq:pre_reduced_Hamiltonian_model}
    A\dot{z} = \mathbb{J}_{2n}\nabla_y H (Az) + r,
\end{equation}
with $r$ being the residual term.
Utilizing the chain rule and the second property of $A^{+}$ in Proposition \ref{prop:properties_symplectic_inverse}, the gradient of the Hamiltonian in \eqref{eq:pre_reduced_Hamiltonian_model} can be recast as
\begin{equation*}
    \nabla_y H(Az) = (A^{+})^{\top} \nabla_z H_{RB} (z).
\end{equation*}
By assuming that the projection residual is orthogonal with respect to the symplectic bilinear form to the space spanned by $A$, we recover
\begin{equation}\label{eq:reduced_Hamiltonian_model}
    \begin{cases}
    \dot{z}(t) = \mathbb{J}_{2k}\nabla_{z}H_{RB}(z(t)),\\
    z(0) = A^{+}y_0.
    \end{cases}
\end{equation}
System \eqref{eq:reduced_Hamiltonian_model} is known as a symplectic Galerkin projection of \eqref{eq:hamilton_equation_tbr} onto $\mathcal{A}$. The pre-processing stage consisting of the collection of all the computations required to assemble the basis $A$ is known as the \emph{offline} stage. The numerical solution of the low-dimensional problem \eqref{eq:reduced_Hamiltonian_model} represents the \emph{online} stage, and provides a fast approximation to the solution of the high-fidelity model \eqref{eq:hamilton_equation_tbr} by means of  \eqref{eq:reduced_basis_approximation}. Even though the offline stage is possibly computationally expensive, this splitting is beneficial in a multi-query context, when multiple instances of \eqref{eq:reduced_Hamiltonian_model} have to be solved, e.g. for parametric PDEs.\\
Traditional projection-based reduction techniques do not guarantee stability, even if the high-dimensional problem admits a stable solution \cite{stability_problem_2}, often resulting in a blowup of system energy. On the contrary, by preserving the geometric structure of the problem, several stability results hold for the reduced Hamiltonian equation \eqref{eq:reduced_Hamiltonian_model}. In \cite[Proposition 15, page A2625]{babak_symplectic}, the authors show that the error in the Hamiltonian $|H(y(t))-H_{RB}(z(t))|$ is constant for all $t$. We detail two relevant results in the following, suggesting that structure and energy preservation are key for stability.
\begin{theorem}[Boundedness result \cite{peng_symplectic}]
Consider the Hamiltonian system \eqref{eq:hamilton_equation_tbr}, with Hamiltonian $H\in C^{\infty}(\mathcal{M})$ and  initial condition $y_0\in\mathbb{R}^{2n}$ such that $y_0\in\text{range}(A)$, with $A\in\mathbb{R}^{2n\times 2k}$ symplectic basis. Let \eqref{eq:reduced_Hamiltonian_model} be the reduced Hamiltonian system obtained as the symplectic Galerkin projection induced by $A$ of \eqref{eq:hamilton_equation_tbr}. If there exists a bounded neighborhood $\mathcal{U}_{y_0}$ in $\mathbb{R}^{2n}$ such that $H(y_0)<H(\Tilde{y})$, or $H(y_0)>H(\Tilde{y})$, for all $\Tilde{y}$ on the boundary of $\mathcal{U}_{y_0}$, then both the original system and the reduced system constructed by the symplectic projection are uniformly bounded for all $t$.
\end{theorem}
\begin{theorem}[Lyapunov stability \cite{peng_symplectic,babak_symplectic}]
Consider the Hamiltonian system \eqref{eq:hamilton_equation_tbr} with Hamiltonian $H\in C^{2}(\mathcal{M})$ and the reduced Hamiltonian system
\eqref{eq:reduced_Hamiltonian_model}. Suppose that $y^{*}$ is a strict local minimum of $H$. Let $S$ be an open ball around $y^{*}$ such that $\nabla^{2}H(y)>0$ and $H(z)<c$, for all $z\in S$ and some $c\in\mathbb{R}$,  and $H(\bar{y})=c$ for some $\bar{y}\in\partial S$, where $\partial S$ is the boundary of $S$. If there exists an open neighborhood $S$ of $y^{*}$ such that $S\cap \text{range}(A)\neq \emptyset$, then the reduced system \eqref{eq:reduced_Hamiltonian_model} has a stable equilibrium point in $S\cap \text{range}(A)$.
\end{theorem}
For the time-discretization of \eqref{eq:reduced_Hamiltonian_model}, the use of a symplectic integrator \cite{hairer_geometric} is crucial for preserving the symplectic structure at the discrete level. In particular, the discrete flow obtained using a symplectic  integrator satisfies a discrete version of Proposition \ref{prop:flow_symplectic_transform}.\\
In the next section, we introduce different strategies to construct symplectic bases as results of optimization problems.

\section{Proper symplectic decomposition}\label{sec:PSD}
Let us consider the solution vectors $y_i=y(t_i)\in\mathbb{R}^{2n}$ (the so-called solution snapshots) obtained, for different time instances $t_i\in [t_0,t_{\text{end}}]$, $\forall i=1,\dots,N$, by time discretization of \eqref{eq:hamilton_equation_tbr} using a symplectic integrator. Define the snapshot matrix
\begin{equation}\label{eq:original_snapshot_matrix}
    M_{y}:=[\enspace y_1 \enspace \dots \enspace y_N \enspace ],
\end{equation}
as the matrix collecting the solution snapshots as columns. In the following, we consider different algorithms stemming from the historical \emph{method of snapshots} \cite{method_of_snapshots}, as the base of the proper orthogonal decomposition (POD).
To preserve the geometric structure of the original model, we focus on a similar optimization problem, the proper symplectic decomposition (PSD), which represents a data-driven basis generation procedure to extract a symplectic basis from $M_y$. It is based on the minimization of the projection error of $M_y$ on $\mathcal{A}$ and it results in the following optimization problem for the definition of the symplectic basis $A\in\mathbb{R}^{2n\times 2k}$:
\begin{equation}\label{eq:symplectic_minimization}
\begin{aligned}
    &\underset{A\in\mathbb{R}^{2n\times 2k}}{\text{minimize }}\quad\| M_y - AA^{+}M_y\|_F, \\
    &\text{subject to }\quad A\in Sp(2k,\mathbb{R}^{2n}),
\end{aligned}
\end{equation}
with $\mathcal{A}=\text{range}(A)$ and $\| \cdot \|_F$ is the Frobenius norm.
Problem \eqref{eq:symplectic_minimization} is similar to the POD minimization, but with the feasibility set of rectangular orthogonal matrices, also known as the Stiefel manifold
\begin{equation*}
\text{St}(2k,\mathbb{R}^{2n}):=\{ L\in\mathbb{R}^{2n\times 2k}: L^{\top}L = \mathbb{I}_{2n}\},
\end{equation*}
replaced by the symplectic Stiefel manifold. Recently there has been a great interest in optimization on symplectic manifolds, and a vast literature is available on the minimization of the least-squares distance from optimal symplectic Stiefel manifolds. This problem has relevant implications in different physical applications, such as the study of optical systems \cite{optical_system} and the optimal control of quantum symplectic gates \cite{optimal_control_gates}.
Unfortunately, with respect to POD minimization, problem \eqref{eq:symplectic_minimization} is significantly more challenging for different reasons. The non-convexity of the feasibility set and the unboundedness of the solution norm precludes standard optimization techniques. Moreover, most of the attention is focused on the case $n=k$, which is not compatible with the reduction goal of MOR.\\
Despite the interest in the topic, an efficient optimal solution algorithm has yet to be found for the PSD. 
Suboptimal solutions have been attained by focusing on the subset of the ortho-symplectic matrices, i.e.,
\begin{equation}
\mathbb{S}(2k,2n):= St(2k,\mathbb{R}^{2n}) \cap Sp(2k,\mathbb{R}^{2n}).
\end{equation}
In \cite{peng_symplectic}, while enforcing the  additional orthogonality constraint 
in \eqref{eq:symplectic_minimization}, the optimization problem is further simplified by assuming a specific structure for $A$. An efficient greedy method, not requiring any additional block structures to $A$, but only its orthogonality and simplecticity, has been introduced in \cite{babak_symplectic}. More recently, in \cite{haasdonk_nonorthonormal_SVD}, the orthogonality requirement has been removed, and different solution methods to the PSD problem are explored. In the following, we briefly review the abovementioned approaches.
\subsection{SVD-based methods for orthonormal symplectic basis generation}\label{sec:svd_based_basis_generators}
In \cite{peng_symplectic}, several algorithms have been proposed to directly construct ortho-symplectic bases.
Exploiting the SVD decomposition of rearranged snapshots matrices, the idea is to search for optimal matrices in subsets of $Sp(2k,\mathbb{R}^{2n})$. Consider the more restrictive feasibility set
\begin{equation*}
    \mathbb{S}_{1}(2k,2n):= Sp(2k,\mathbb{R}^{2n}) \cap 
    \Bigg \{
    \begin{bmatrix}
    \enspace \Phi &\enspace 0 \enspace\\
    \enspace 0 &\enspace \Phi \enspace 
    \end{bmatrix} \Bigg |
    \; \Phi\in\mathbb{R}^{n\times k}
    \Bigg \}.
\end{equation*}
Then $A^{\top}\mathbb{J}_{2n}A = \mathbb{J}_{2k}$ holds if and only if $\Phi^{\top}\Phi=\mathbb{I}_n$, i.e. $\Phi\in St(k,\mathbb{R}^{n})$. Moreover, we have that $A^{+}=\text{diag}(\Phi^{\top},\Phi^{\top})$. The cost function in \eqref{eq:symplectic_minimization} becomes
\begin{equation}\label{eq:symplectic_minimization_CL}
\| M_y - AA^{+}M_y\|_F = \| M_1 - \Phi \Phi^{\top}M_1\|_F,
\end{equation}
with $M_1=[ \enspace p_1 \enspace \dots \enspace p_N \enspace q_1 \enspace \dots \enspace q_N \enspace ]\in\mathbb{R}^{n\times 2N}$, where $p_i$ and $q_i$ are the generalized phase-space components of $y_i$. Thus, as a result of the Eckart–Young–Mirsky theorem,  \eqref{eq:symplectic_minimization_CL} admits a solution in terms of the singular-value decomposition of the data matrix $M_1$. 
This algorithm, formally known as Cotangent Lift, owes its name to the interpretation of the solution $A$ to \eqref{eq:symplectic_minimization_CL} in $\mathbb{S}_1(2k,2n)$ as the cotangent lift of linear mappings, represented by $\Phi$ and $\Phi^{\top}$, between vector spaces of dimensions $n$ and $k$. Moreover, this approach constitutes the natural outlet in the field of Hamiltonian systems of the preliminary work of Lall et al. \cite{lall_lagrangian} on structure-preserving reduction of Lagrangian systems.
However, there is no guarantee that the Cotangent Lift basis is close to the optimal of the original PSD functional. \\
A different strategy, known as Complex SVD decomposition, relies on the definition of the complex snapshot matrix $M_2 = [\enspace p_1+iq_1 \enspace \dots \enspace p_N+iq_N\enspace ]\in\mathbb{C}^{n\times N}$, with $i$ being the imaginary unit. Let $U=\Phi+i\Psi \in\mathbb{C}^{n\times N}$, with $\Phi,\Psi\in\mathbb{R}^{n\times k}$, be the unitary matrix solution to the following accessory problem:
\begin{equation}\label{eq:minimization_complex_SVD}
\begin{aligned}
    &\underset{U\in\mathbb{R}^{n\times 2k}}{\text{minimize }}\quad\| M_2 - UU^{*}M_2\|_F, \\
    &\text{subject to }\quad U\in St(2k,\mathbb{R}^{2n}).
\end{aligned}
\end{equation}
As for the Cotangent Lift algorithm, the solution to \eqref{eq:minimization_complex_SVD} is known to be the set of the $2k$ left-singular vectors of $M_2$ corresponding to its largest singular values. In terms of the real and imaginary parts of $U$, the orthogonality constraint implies
\begin{equation}\label{eq:conditions_complex_SVD}
    \Phi^{\top}\Phi + \Psi^{\top}\Psi = \mathbb{I}_n, \qquad \Phi^{\top}\Psi = \Psi^{\top}\Phi.
\end{equation}
Consider the ortho-symplectic matrix, introduced in \cite{peng_symplectic}, and given by
\begin{equation}\label{eq:general_orthosymplectic}
    A = 
    \begin{bmatrix}
    \enspace E &\enspace \mathbb{J}_{2n}^{\top}E \enspace
    \end{bmatrix}
    \in\mathbb{R}^{2n\times 2k},
    \qquad
    E^{\top}E = \mathbb{I}_k, \quad E^{\top}\mathbb{J}_{2n}E = 0_{k},
    \qquad \text{with} \quad
    E = \begin{bmatrix}
    \enspace \Phi \enspace \\
    \enspace \Psi \enspace
    \end{bmatrix}.\end{equation}
Using \eqref{eq:conditions_complex_SVD}, it can be shown that such an $A$ is the optimal solution of the PSD problem in
\begin{equation*}
    \mathbb{S}_{2}(2k,2n):= Sp(2k,\mathbb{R}^{2n}) \cap 
    \Bigg \{
    \begin{bmatrix}
    \enspace \Phi \;&\enspace -\Psi \enspace\\
    \enspace \Psi &\enspace \Phi \enspace 
    \end{bmatrix} \Bigg |
    \; \Phi,\Psi\in\mathbb{R}^{n\times k}
    \Bigg \},
\end{equation*}
 that minimizes the projection error of $M_r := [\enspace M_y \enspace \mathbb{J}_{2n}M_y \enspace]$, also known as the rotated snapshot matrix, with $M_y$ given in \eqref{eq:original_snapshot_matrix}. In \cite{haasdonk_nonorthonormal_SVD}, extending the result obtained in \cite{symplectic_schur} for square matrices, it has been shown that  \eqref{eq:general_orthosymplectic} is a complete characterization of symplectic matrices with orthogonal columns, meaning that all the ortho-symplectic matrices admit a representation of the form \eqref{eq:general_orthosymplectic}, for a given $E$, and hence $\mathbb{S}_2(2k,2n) \equiv \mathbb{S}(2k,2n)$. In the same work, Haasdonk et al. showed that an ortho-symplectic matrix that minimizes the projection error of $M_r$ is also a minimizer of the projection error of the original snapshot matrix $M_y$, and vice versa. This is been achieved by using an equivalence argument based on the POD applied to the matrix $M_r$. Thus, combining these two results, the Complex SVD algorithm provides a minimizer of the PSD problem for ortho-symplectic matrices.
\subsection{SVD-based methods for non-orthonormal symplectic basis generation}
In the previous section, we showed that the basis provided by the Complex SVD method is not only near-optimal in $\mathbb{S}_2$, but is optimal for the cost functionals in the space of ortho-symplectic matrices. 
The orthogonality of the resulting basis is beneficial  \cite{Kressner_condition}, among others, for reducing the condition number associated with the fully discrete formulation of \eqref{eq:reduced_Hamiltonian_model}. A suboptimal solution to the PSD problem not requiring the orthogonality of the feasibility set is proposed in \cite{peng_symplectic}, as an improvement of the SVD-based generators of ortho-symplectic bases using the Gappy POD \cite{gappy_POD}, under the name of nonlinear programming approach (NLP). Let $A^*\in\mathbb{S}_2(2r,2n)$ be a basis  of dimension $2r$ generated using the Complex SVD method. The idea of the NLP is to construct a target basis $A\in Sp(2k,\mathbb{R}^{2n})$, with $k<r \ll n$, via the linear mapping
\begin{equation}\label{eq:nonlinear_programming_map}
A = A^{*}C,
\end{equation}
with $C\in\mathbb{R}^{2r\times 2k}$. Using \eqref{eq:nonlinear_programming_map} in \eqref{eq:symplectic_minimization} results in a PSD optimization problem for the coefficient matrix $C$, of significant smaller dimension ($4kr$ parameters) as compared to the original PSD problem ($4kn$ parameters) with $A$ unknown. However, no optimality results are available for the NLP method.\\
A different direction has been pursued in \cite{haasdonk_nonorthonormal_SVD}, based on the connection between traditional SVD and Schur forms and the matrix decompositions, related to symplectic matrices, as proposed in the following theorem.  
\begin{theorem}[SVD-like decomposition {\cite[Theorem 1, page 6]{svd_like_decomposition}}]
If $B\in\mathbb{R}^{2n\times n_s}$, then there exists $S\in Sp(2n,\mathbb{R}^{2n})$, $Q\in St(n_s,\mathbb{R}^{n_s})$ and $D\in\mathbb{R}^{2n\times n_s}$ of the form
\begin{equation}\label{eq:D_matrix}
D = 
    \begin{blockarray}{ccccc}
     \enspace b & q & b & n-2b-q & \\
      \begin{block}{[cccc]c}
      \Sigma & 0 & 0 & 0 & b \\ 
      0 & I & 0 & 0 & q \\
      0 & 0 & 0 & 0 & m-b-q \\
      0 & 0 & \Sigma & 0 & b \\
      0 & 0 & 0 & 0 & q \\
      0 & 0 & 0 & 0 & m-b-q \\
      \end{block}
    \end{blockarray},
\end{equation}
with $\Sigma=\text{diag}(\sigma_1,\dots,\sigma_b)$, $\sigma_i>0$ $\forall i=1,\dots,b$, such that
\begin{equation}\label{eq:symplectic_decomposition}
B=SDQ.
\end{equation}
Moreover, rank$(B)=2b+q$ and $\sigma_i$ are known as symplectical singular values.
\end{theorem}
\noindent Let us apply the SVD-like decomposition to the snapshot matrix $M_y$ \eqref{eq:original_snapshot_matrix}, where $n_s$ represents the number of snapshots, and define its weighted symplectic singular values as 
\begin{equation*}
w_i =
\begin{cases}
\sigma_i \sqrt{\| S_i \|_2^2+\| S_{n+i}\|_2^2 },  & 1\leq i \leq b,\\
\| S_i \|_2, & b+1\leq i \leq b+q,
\end{cases}
\end{equation*} 
with $S_i\in\mathbb{R}^{2n}$ being the $i$-th column of $S$ and $\| \cdot \|_2$ the Euclidean norm.
The physical interpretation of the classical POD approach characterizes the POD reduced basis as the set of a given cardinality that captures most of the energy of the system. The energy retained in the reduced approximation is quantified as the sum of the squared singular values corresponding to the left singular vectors of the snapshot matrix representing the columns of the basis.
A similar guiding principle is used in \cite{haasdonk_nonorthonormal_SVD}, where the energy of the system, i.e., the Frobenius norm of the snapshot matrix, is connected to the weighted symplectic singular values as
\begin{equation}\label{eq:energy_symplectic_values}
\| M_y \|_F^2 = \sum_{i=1}^{b+q} w_i^2.
\end{equation}
 Let $\mathcal{I}_{\text{PSD}}$ be the set of indices corresponding to the $k$ largest energy contributors in  \eqref{eq:energy_symplectic_values},
\begin{equation}
\mathcal{I}_{\text{PSD}} = \{ i_j \}_{j=1}^{k} = \underset{\mathcal{I}\subset \{ 1,\dots,b+q\}}{\text{argmax}}\Big(\sum_{i\in\mathcal{I}} w_i^2 \Big).
\end{equation}
Then, the PSD SVD-like decomposition defines a symplectic reduced basis $A\in Sp(2k,\mathbb{R}^{2n})$ by selecting the pairs of columns from the symplectic matrix $S$ corresponding to the indices set $\mathcal{I}_{PSD}$
\begin{equation}
A=[\enspace s_{i_1} \enspace \dots\enspace s_{i_k} \enspace s_{n+i_1} \enspace \dots \enspace s_{n+i_k} \enspace].
\end{equation}
Similarly to the POD, the reconstruction error of the snapshot matrix depends on the magnitude of the discarded weighted symplectic singular values as
\begin{equation}
\| M_y - AA^{+}M_y \|_F^{2} = \sum_{i\in\{1,\dots,b+q\} \setminus\mathcal{I}_{\text{PSD}}}w_i^{2}.
\end{equation}
Even though there are no proofs that the PSD SVD-like algorithm reaches the global optimum in the sense of \eqref{eq:symplectic_minimization}, some  analysis and numerical investigations  suggest that it provides superior results as compared to orthonormal techniques \cite{haasdonk_nonorthonormal_SVD}.

\subsection{Greedy approach to symplectic basis generation}\label{sec:greedy_approach}
The reduced basis methodology is motivated and applied within the context of real-time and multi-queries simulations of parametrized PDEs. In the framework of Hamiltonian systems, we consider the following parametric form of \eqref{eq:hamilton_equation_tbr}
\begin{equation}\label{eq:param_hamilton_equation_tbr}
    \begin{cases}
    \dot{y}(t,\mu)=\mathbb{J}_{2n}\nabla_y H(y(t,\mu);\mu),\\
    y(0,\mu) = y_{0}(\mu),
    \end{cases}
\end{equation}
with $\mu\in\mathcal{P}\subset \mathbb{R}^{d}$ being a $d$-dimensional parameter space. Let $\mathcal{Z}^{\mathcal{P}}$ be the set of solutions to \eqref{eq:param_hamilton_equation_tbr} defined as
\begin{equation*}
\mathcal{Z}^{\mathcal{P}} = \{ y(t,\mu):t\in [t_0,t_{\text{end}}], \mu\in\mathcal{P} \}\subset \mathbb{R}^{2n}.
\end{equation*}
For the sake of simplicity, in the previous sections we have only considered the non-parametric case. The extension of SVD-based methods for basis generations to \eqref{eq:param_hamilton_equation_tbr} is straightforward on paper, but it is often computationally problematic in practice as the number of snapshots increases. Similar to other SVD-based algorithms, the methods described in the previous sections require the computation of the solution to \eqref{eq:param_hamilton_equation_tbr} corresponding to a properly chosen discrete set of parameters $S^{\mu}=\{ \mu_j\}_{j=1}^{p}\subset \mathcal{P}$ and time instances $S^{t}=\{ t_i \}_{i=1}^{N}$, defined \emph{a priori},
and constituting the sampling set $S^{\mu,t}:=S^{\mu}\times S^{t}$. Random or structured strategies exist to define the set $S^{\mu}$, such as the Monte Carlo sampling, Latin hypercube sampling, and sparse grids \cite{sampling_techniques_1}, while $S^{t}$ is a subset of the time-discretization, usually dictated by the integrator of choice. The set of snapshots corresponding to the sampling set $S^{\mu,t}$ must provide a ``good'' approximation of the solution manifold and should not miss relevant parts of the time-parameter domain.
Once the sampling set $S^{\mu,t}$ has been fixed, the matrix
$M_y$, $M_1$, or $M_2$, depending on the method of choice, is assembled, and its singular value decomposition is computed. Even though a certain amount of computational complexity is tolerated in the \emph{offline} stage to obtain a significant speed-up in the \emph{online} stage, the evaluation of the high-fidelity solution for a large sampling set and the SVD of the corresponding snapshot matrix are often impractical or not even feasible. 
Hence, an efficient approach is an incremental procedure. The reduced basis, in which the column space represents the approximating manifold, is improved iteratively by adding basis vectors as columns. The candidate basis vector is chosen as the maximizer of a much cheaper optimization problem. This summarizes the philosophy of the greedy strategy applied to RB methods
\cite{strong_greedy_exponential_convergence,strong_greedy_algebraic_convergence}, which requires two main ingredients: the definition of an error indicator and a process to add a candidate column vector to the basis.\\
Let $U_k$ be an orthonormal reduced basis produced after $k$ steps of the algorithm. In its idealized form, introduced in \cite{standard_greedy}, the greedy algorithm uses the projection error
\begin{equation}\label{eq:error_greedy}
    (t^{*},\mu^{*}) := \underset{(t_i,\mu_{j})\in S^{\mu,t}}{\text{argmax}} \| u(t_i,\mu_{j}) - U_k U_k^{\top} u(t_i,\mu_{j}) \|_2,
\end{equation}
 to identify the snapshot $u^{*}:=u(t^{*},\mu^{*})$ that is worst approximated by the column space of $U_k$ over the entire sampling set $S^{\mu,t}$. Let $u_{k+1}$ be the vector obtained by orthonormalizing $u^{*}$ with respect to $U_k$. Then the basis $U_k$ is updated as $U_{k+1}=[ \enspace U_k \enspace u_{k+1} \enspace]$.  To avoid the accumulation of rounding errors, it is preferable to utilize backward stable orthogonalization processes, such as the modified Gram-Schmidt orthogonalization. The algorithm terminates when the basis reaches the desired dimension, or the error \eqref{eq:error_greedy} is below a certain tolerance. In this sense, the basis $U_{k+1}$ is \emph{hierarchical} because its column space contains the column space of its previous iterations. This process is referred to as \emph{strong greedy} method. Even though introduced as a heuristic procedure, interesting results regarding algebraic and exponential convergence have been formulated in \cite{strong_greedy_exponential_convergence,strong_greedy_algebraic_convergence}, requiring the orthogonality of the basis in the corresponding proofs. 
However, in this form, the scheme cannot be efficiently implemented: the error indicator  \eqref{eq:error_greedy} is expensive to calculate because it requires all the snapshots of the training set $S^{\mu,t}$ to be accessible, relieving the computation only of the cost required for the SVD.\\
An adjustment of the \emph{strong greedy} algorithm, known as \emph{weak greedy} algorithm, assembles  the snapshot matrix corresponding to $S^{\mu,t}$ iteratively while expanding the approximating basis. The idea is to replace \eqref{eq:error_greedy} with a surrogate indicator $\eta:S^{\mu,t}\mapsto \mathbb{R}$ that does not demand the computation of the high-fidelity solution for the entire time-parameter domain.\\
In the case of elliptic PDEs, an \emph{a-posteriori} residual-based error indicator requiring a polynomial computational cost in the approximation space dimension has been introduced in \cite{residual_based_indicator}. The substantial computational savings allow the choice of a more refined, and therefore representative, sampling set $S^{\mu,t}$. One might also use a goal-oriented indicator as the driving selection in the greedy process to obtain similar computational benefits. In this direction, in the framework of structure-preserving model order reduction, \cite{babak_symplectic} suggests the Hamiltonian as a proxy error indicator. Suppose $A_{2k}=[ \enspace E_k \enspace\mathbb{J}_{2n}^{\top}E_k \enspace]$, with $E_k = [ \enspace e_1 \enspace \dots \enspace e_k \enspace]$, is a given ortho-symplectic basis and consider
\begin{equation}\label{eq:Hamiltonian_indicator_greedy}
    (t^{*},\mu^{*}) := \underset{(t_i,\mu_{j})\in S^{\mu,t}}{\text{argmax}} | H(y(t_i,\mu_{j})) - H(A_{2k} A_{2k}^{+} y(t_i,\mu_{j})) |.
    \end{equation}
By \cite[Proposition 15]{babak_symplectic}, the error in the Hamiltonian depends only on the initial condition and the symplectic reduced basis. Hence, the indicator \eqref{eq:Hamiltonian_indicator_greedy} does not require integrating in time the full system \eqref{eq:param_hamilton_equation_tbr} over the entire set $S^{\mu}$, but only over a small fraction of the parameter set, making the procedure fast. Hence, the parameter space can be explored first,
\begin{equation}
\mu^{*} := \underset{\mu_j\in S^{\mu}}{\text{argmax}}| H(y_0(\mu_j))-H(A_{2k}A_{2k}^{+}y_0(\mu_j)) |,
\end{equation}
to identify the value of the parameter that maximizes the error in the Hamiltonian as a function of the initial condition. This step may fail if $y_{0}(\mu_j)\in\text{range}(A_{2k})$, $\forall j=1,\dots,p$. Then \eqref{eq:param_hamilton_equation_tbr} is temporally integrated to collect the snapshot matrix
$M_g = [\enspace y(t_1,\mu^{*}) \enspace \dots \enspace  y(t_N,\mu^{*}) \enspace ]$.
Finally, the candidate basis vector $y^{*}=y(\mu^{*},t^{*})$ is selected as the snapshot that maximizes the projection error 
\begin{equation}
t^{*}:=\underset{t_i\in S^{t}}{\text{argmax}} \| y(t_i,\mu^{*})-A_{2k}A_{2k}^{+}y(t_i,\mu^{*})\|_2.
\end{equation}
Standard orthogonalization techniques, such as QR methods, fail to preserve the symplectic structure \cite{instability_QR}. In \cite{babak_symplectic}, the SR method \cite{SR}, based on the symplectic Gram-Schimidt, is employed to compute the additional basis vector $e_{k+1}$ that conforms to the geometric structure of the problem. To conclude the $(k+1)$-th iteration of the algorithm, the basis $A_{2k}$ is expanded in 
\begin{equation*}
A_{2(k+1)} = [ \enspace E_{k} \enspace e_{k+1} \enspace \mathbb{J}_{2n}^{\top}E_k \enspace \mathbb{J}_{2n}^{\top}e_{k+1} \enspace ].
\end{equation*}
We stress that, with this method, known as symplectic greedy RB, two vectors, $e_{k+1}$ and $\mathbb{J}_{2n}^{\top}e_{k+1}$, are added to the symplectic basis at each iteration, because of the structure of ortho-symplectic matrices. A different strategy, known as PSD-Greedy algorithm and partially based on the PSD SVD-like decomposition, has been introduced in \cite{haasdonk_greedy}, with the feature of not using orthogonal techniques to compress the matrix $M_g$.
In \cite{babak_symplectic}, following the results given in \cite{strong_greedy_exponential_convergence}, the exponential convergence of the symplectic strong greedy method has been proved.
\begin{theorem}[{\cite[Theorem 20, page A2632]{strong_greedy_exponential_convergence}}]\label{th:convergence_symplectic_greedy}
Let $\mathcal{Z}^{\mathcal{P}}$ be a compact subset of $\mathbb{R}^{2n}$. Assume that the Kolmogorov $m$-width of $\mathcal{Z}^{\mathcal{P}}$ defined as
\begin{equation*}
d_m(\mathcal{Z}^{\mathcal{P}}) = 
\underset{\substack{\mathcal{Z}_*\subset\mathbb{R}^{2n}\\ \text{dim}(\mathcal{Z}_{*})=m}}{\text{inf}}\;
\underset{v\in\mathcal{Z}^{\mathcal{P}}}{\text{sup}} \;
\underset{w\in\mathcal{Z_*}}{\text{min}}
\| v-w \|_2, 
\end{equation*}
decays exponentially fast, namely $d_m(\mathcal{Z}^{\mathcal{P}})\leq c \exp (-\alpha m)$ with $\alpha> \log 3$. Then there exists $\beta>0$ such that the symplectic basis $A_{2k}$ generated by the symplectic strong greedy algorithm provides exponential approximation properties,
\begin{equation}\label{eq:convergenge_symplectic_greedy_method}
\| s - A_{2k}A_{2k}^{+} s\|_2 \leq C \exp(-\beta k),
\end{equation}
for all $s\in \mathcal{Z}^{\mathcal{P}}$ and some $C>0$.
\end{theorem}
Theorem \ref{th:convergence_symplectic_greedy} holds only when the projection error is used as the error indicator instead of the error in the Hamiltonian. However, it has been observed for different symplectic parametric problems
\cite{babak_symplectic} that the symplectic method using the loss in the Hamiltonian converges  with the same rate of \eqref{eq:convergenge_symplectic_greedy_method}.  The orthogonality of the basis is used to prove the convergence of the greedy procedure. In the case of a non-orthonormal symplectic basis, supplementary assumptions are required to  ensure the convergence of the algorithm.

\section{Dynamical low-rank reduced basis methods for Hamiltonian systems}\label{sec:dynamical_low_rank}
The Kolmogorov m-width of a compact set describes how well this can be approximated by a linear subspace of a fixed dimension $m$. A problem \eqref{eq:param_hamilton_equation_tbr} is informally defined \emph{reducible} if $d_m$ decays sharply with $m$, implying the existence of a low-dimensional representation of $\mathcal{Z}^{\mathcal{P}}$. A slow decay limits the accuracy of any efficient projection-based reduction on linear subspaces, including all the methods discussed so far. For Hamiltonian problems, often characterized by the absence of physical dissipation due to the conservation of the Hamiltonian, we may have $d_m(\mathcal{Z}^{\mathcal{P}})=\mathcal{O}(m^{-\frac{1}{2}})$ in case of discontinuous initial condition \cite{decay_wave} for wave-like problems.
Several techniques, either based on nonlinear transformations of the solution manifold to a reducible framework \cite{ohlberger_nonlinear_manifold} or presented as online adaptive methods to target solution manifolds at fixed time \cite{willcox_update}, have been introduced to overcome the limitations of the linear approximating spaces. In different ways, they all abandon the framework of symplectic vector spaces. Therefore, none of them guarantees conservation of the symplectic structure in the reduction process. Musharbash et al. \cite{nobile_dlr} proposed a dynamically orthogonal (DO) discretization of stochastic wave PDEs with a symplectic structure. In the following, we outline the structure-preserving dynamic RB method for parametric Hamiltonian systems, proposed by Pagliantini \cite{pagliantini} in the spirit of the geometric reduction introduced in \cite{feppon}. In contrast with traditional methods that provide a global basis, which is fixed in time, the gist of a dynamic approach is to evolve a local-in-time basis to provide an accurate approximation of the solution to the parametric problem \eqref{eq:param_hamilton_equation_tbr}. The idea is to exploit the  local low-rank nature of Hamiltonian dynamics in the parameter space.
From a geometric perspective, the approximate solution evolves according to naturally constrained dynamics, rather than weakly enforcing the required properties, such as orthogonality or symplecticity of the RB representation, via Lagrange multipliers. This result is achieved by viewing the flow of the reduced model as prescribed by a vector field that is everywhere tangent to the desired manifold. \\
Suppose we are interested in solving \eqref{eq:param_hamilton_equation_tbr} for a set of $p$ vector-valued parameters $\eta_h=\{\mu_i\}_{i=1}^{p}$, sampled from $\mathcal{P}$. Then the Hamiltonian system, evaluated at $\eta_h$, can be recast as a set of ODEs in the matrix unknown $\mathcal{R}\in\mathbb{R}^{2n\times p}$,
\begin{equation}\label{eq:system_Hamiltonian}
\begin{cases}
\dot{\mathcal{R}}(t) = \mathcal{X}_H(\mathcal{R}(t),\eta_h) = \mathbb{J}_{2n} \nabla_{\mathcal{R}} H(\mathcal{R}(t),\eta_h),\\
\mathcal{R}(t_0) = \mathcal{R}_0(\mu_h),
\end{cases}
\end{equation}
where $H$ is a vector-valued Hamiltonian function, the $j$-th column of $\mathcal{R}(t)$ is such that $\mathcal{R}_j(t) = y(t,\mu_j)$,
and $(\nabla_{\mathcal{R}} H)_{i,j}:=\partial H_j/ \partial \mathcal{R}_{i,j}$.
We consider an approximation of the solution to \eqref{eq:system_Hamiltonian} of the form
\begin{equation}\label{eq:dynamical_reduced_basis}
\mathcal{R}(t)\approx R(t)=A(t)Z(t),
\end{equation}
where $A(t)\in\mathbb{S}(2k,2n)$, and $Z(t)\in\mathbb{R}^{2k\times p}$ is such that its $j$-th column $Z_j(t)$ collects coefficients, with respect to the basis $A(t)$, of the approximation of $y(t,\mu_j)$. Despite being cast in the same framework of an RB approach, a stark difference between \eqref{eq:dynamical_reduced_basis} and \eqref{eq:reduced_basis_approximation} lies in the time-dependency of the basis in \eqref{eq:dynamical_reduced_basis}.\\
Consider the manifold of $2n\times p$ matrices having at most rank $2k$, and defined as
\begin{equation}
\mathcal{Z}^{\mathcal{P}}_{2n}:=\{ R\in\mathbb{R}^{2n\times p}: R=AZ \; \text{ with } \; A\in\mathbb{S}(2k,2n), \; Z\in\mathbb{Z}  \},
\end{equation}
with the technical requirement 
\begin{equation}
\mathbb{Z}:=\{ Z\in\mathbb{R}^{2k\times p}: \text{rank}(ZZ^{\top}+\mathbb{J}_{2k}^{\top}ZZ^{\top}\mathbb{J}_{2k}) = 2k\}.
\end{equation}
This represents a full-rank condition on $Z$ to ensure uniqueness of the representation \eqref{eq:dynamical_reduced_basis} for a fixed basis.
The tangent vector at $R(t)=A(t)Z(t)\in\mathcal{Z}^{\mathcal{P}}_{2n}$ is given by $X=X_{A}Z+AX_{Z}$, where $X_A$ and $X_Z$ correspond to the tangent directions for the time-dependent matrices $A$ and $Z$, respectively. Applying the orthogonality and symplecticity condition on $A(t)$,
%$A^{\top}A=I_{2n}$ and $A^{\top}\mathbb{J}_{2n}A=\mathbb{J}_{2k}$ must hold
for all times $t$, results in 
\begin{equation}\label{eq:time_constraint}
X_A^{\top}A+A^{\top}X_A=0 \enspace \text{and} \enspace X_A^{\top}\mathbb{J}_{2n}A + A^{\top}\mathbb{J}_{2n}X_A=0,
\end{equation}
respectively. Using \eqref{eq:time_constraint} and an additional gauge constraint to uniquely parametrize the tangent vectors $X$ by the displacements $X_A$ and $X_Z$, the tangent space of $\mathcal{Z}^{\mathcal{P}}_{2n}$ at $R=AZ$ can be characterized as 
\begin{equation*}
\begin{aligned}
T_{R}\mathcal{M}^{\mathcal{P}}_{2n} = \{ X\in\mathbb{R}^{2n\times p}:  &X=X_AZ + AX_Z,\\
 &\text{with }X_Z\in\mathbb{R}^{2k\times p}, X_A\in\mathbb{R}^{2n\times 2k}, X_A^{\top}A = 0, X_A\mathbb{J}_{2k} = \mathbb{J}_{2n}X_A \}.
\end{aligned}
\end{equation*}
The reduced flow describing the evolution of the approximation $R(t)$ is derived in \cite{pagliantini} by projecting the full velocity field $\mathcal{X}_H$ in \eqref{eq:system_Hamiltonian} onto the tanget space $T_{R(t)}\mathcal{Z}^{\mathcal{P}}_{2n}$ of $\mathcal{Z}^{\mathcal{P}}_{2n}$ at $R(t)$, i.e.,
\begin{equation}\label{eq:symplectic_projection_dynamics}
\begin{cases}
\dot{R}(t)=\Pi_{T_{R(t)}\mathcal{Z}^{\mathcal{P}}_{2n}} \mathcal{X}_H(R(t),\eta_h),\\
R(t_0) = U_0Z_0.
\end{cases}
\end{equation}
To preserve the geometric structure of the problem, the projection operator $\Pi_{T_{R(t)}\mathcal{Z}^{\mathcal{P}}_{2n}}$  is a symplectomorphism (see Definition \ref{def:symplectic_transform}) for each realization of the parameter $\mu_j\in \eta_h$, in the sense given in the following proposition.
\begin{proposition}[{\cite[Proposition 4.3, page 420]{pagliantini}}]
Let $S:=ZZ^{\top}+\mathbb{J}_{2k}ZZ^{\top}\mathbb{J}_{2k}\in\mathbb{R}^{2k\times 2k}$. Then, the map
\begin{equation}\label{eq:explicit_definition_projection}
\begin{aligned}
\Pi_{T_{R(t)}\mathcal{Z}^{\mathcal{P}}_{2n}}: \enspace &\mathbb{R}^{2n\times p} &&\rightarrow T_{R(t)}\mathcal{Z}^{\mathcal{P}}_{2n},\\
\enspace &w && \mapsto (I_{2n}-AA^{\top})(wZ^{\top}+\mathbb{J}_{2n}wZ^{\top}\mathbb{J}_{2k}^{\top})S^{-1}Z+AA^{\top}w,
\end{aligned}
\end{equation}
is a symplectic projection, in the sense that 
\begin{equation*}
\sum_{j=1}^{p}\Omega^{\mu_j} (w-\Pi_{T_{R(t)}\mathcal{Z}^{\mathcal{P}}_{2n}}w,y) =0, \quad \forall y\in \Pi_{T_{R(t)}\mathcal{Z}^{\mathcal{P}}_{2n}},
\end{equation*}
where $\Omega^{\mu_j}$ is the symplectic form associated with the parameter $\mu_j$.
\end{proposition}
\noindent The optimality of the reduced dynamics, in the Frobenius norm, follows from \eqref{eq:symplectic_projection_dynamics},
where the flow of $R$ is prescribed by the best low-rank approximation of the Hamiltonian velocity field vector $\mathcal{X}_H$ into the tangent space of the reduced manifold $\mathcal{Z}_{2n}^{\mathcal{P}}$.
Using \eqref{eq:explicit_definition_projection} and \eqref{eq:symplectic_projection_dynamics}, it is straightforward to derive the evolution equations for $A(t)$ and $Z(t)$:
\begin{equation}\label{eq:DLR_symplectic}
\begin{cases}
\dot{Z}_j(t) = \mathbb{J}_{2n} \nabla_{Z_j} H(AZ_j,\mu_j),\\
\dot{A}(t)=(I_{2n}-AA^{\top})(\mathbb{J}_{2n}YZ - YZ\mathbb{J}_{2n}^{\top})S^{-1},  \\
A(t_0)Z(t_0) = A_0Z_0,
\end{cases}
\end{equation}
with $Y:=[ \enspace \mathbb{J}_{2n}\nabla H(UZ_1,\mu_1) \enspace \dots \enspace \mathbb{J}_{2n}\nabla H(UZ_p,\mu_p) \enspace ]$.\\
The coefficients $Z$ evolve according to a system of $p$ independent Hamiltonian equations, each in $2n$ unknowns, corresponding to the symplectic Galerkin projection onto $\text{range}(A)$ for each parameter instance in $\eta_h$, similarly to the global symplectic RB method \eqref{eq:reduced_Hamiltonian_model}. In \eqref{eq:DLR_symplectic}, however, the basis $A$ evolves in time according to a matrix equation in $2n\times 2k$ unknowns, affecting the projection.
A crucial property of the structure of $A(t)$ is given in the following proposition.
\begin{proposition}[{\cite[Proposition 4.5, page 423]{pagliantini}}]\label{prop:conservation_DLR}
If $A_0\in\mathbb{S}(2k,2n)$ then $A(t)\in\mathbb{R}^{2n\times 2k}$ solution of \eqref{eq:DLR_symplectic} satisfies $A(t)\in\mathbb{S}(2k,2n)$ for all $t>t_0$.
\end{proposition}
Standard numerical integrators, applied to \eqref{eq:DLR_symplectic}, do not preserve, at the time-discrete level, the property in Proposition \ref{prop:conservation_DLR} and the ortho-symplectic structure is compromised after a single time step. 
In \cite{pagliantini}, two different intrinsic integrators have been investigated to preserve the ortho-symplecticity of the basis, based on Lie groups and tangent techniques. 
Both methods require the introduction of a local chart defined on the tangent space $T_{A(t)}\mathbb{S}$ of the manifold $\mathbb{S}(2k,2n)$ at $A(t)$, with
\begin{equation*}
T_{A(t)}\mathbb{S} := \{ V\in\mathbb{R}^{2n\times 2k}: A^{\top} V\in \mathfrak{g}(2k) \}
\end{equation*}
and $\mathfrak{g}(2k)$ being the vector space of skew-symmetric and Hamiltonian $2k\times 2k$ real square matrices.
In terms of differential manifolds, $\mathfrak{g}(2k)$ represents, together with the Lie bracket $[\cdot,\cdot ]:\mathfrak{g}(2k)\times\mathfrak{g}(2k)\mapsto \mathfrak{g}(2k)$ defined as the matrix commutator
$[M,L]:=ML-LM$, with $M,L\in\mathfrak{g}(2k)$,
the Lie algebra corresponding to the Lie group $\mathbb{S}(2k,2k)$.
The idea is  to recast the basis equation in \eqref{eq:DLR_symplectic} in an evolution equation in the corresponding Lie algebra. The linearity of Lie algebras allows to compute, via explicit Runge-Kutta methods, numerical solutions that remain on the Lie algebra.
Finally, the Cayley transform  $\text{cay}:\mathfrak{g}(2k)\mapsto \mathbb{S}(2k,2k)$ is exploited to generate local coordinate charts and retraction /inverse retraction maps, used to recover the solution in the manifold of rectangular ortho-symplectic matrices. 
In \cite{rank_adaptive_DLR}, the structure-preserving dynamical RB-method has been paired with a rank-adaptive procedure, based on a residual error estimator, to dynamically update also the dimension of the basis. 

\section{Extensions to more general Hamiltonian problems}\label{sec:extension}
\subsection{Dissipative Hamiltonian systems}
Many areas of engineering require a more general framework than the one offered by classical Hamiltonian systems, requiring the inclusion of energy-dissipating elements. 
While the principle of energy conservation is still used to describe the state dynamics, dissipative perturbations must be modelled and introduced in the Hamiltonian formulation \eqref{eq:hamilton_equation_tbr}.
Dissipative Hamiltonian systems, with so-called Rayleigh type dissipation, are considered a special case of forced Hamiltonian systems, with the state $y=(q,p)\in\mathbb{R}^{2n}$, with $q,p\in\mathbb{R}^{n}$, following the time evolution given by
\begin{equation}\label{eq:dissipative_Hamiltonian}
\begin{cases}
\dot{y}(t) = \mathbb{J}_{2n}\nabla H(y(t))+\mathcal{X}_F(y(t)),\\
y(0)=y_0,
\end{cases}
\end{equation}
where $\mathcal{X}_F\in\mathbb{R}^{2n}$ is a velocity field, introducing dissipation, of the form 
\begin{equation}
\mathcal{X}_F := 
\begin{bmatrix}
0_{n}\\
f_H(y(t))
\end{bmatrix}.
\end{equation}
We require $\mathcal{X}_F$ to satisfy $(\nabla_y H)^{\top}\mathcal{X}_F\leq 0$, $\forall y\in\mathbb{R}^{2n}$, to represent a dissipative term and therefore 
\begin{equation}\label{eq:dissipation_condition}
(\nabla_p H)^{\top} f_H\leq 0.
\end{equation}
In terms of Rayleigh dissipation theory, there exists a symmetric positive semidefinite matrix $R(q)\in\mathbb{R}^{n\times n}$ such that $f_H=-R(q)\dot{q}(p,q)$ and \eqref{eq:dissipation_condition} reads 
\begin{equation*}
(\nabla_p H)^{\top} f_H = \dot{q}^{\top}f_H = -\dot{q}^{\top}R(q)\dot{q}\leq0.
\end{equation*}
Several strategies have been proposed to generate stable reduced approximations of \eqref{eq:dissipative_Hamiltonian}, based on Krylov subspaces or POD \cite{port_hamiltonian_reduction_1,port_hamiltonian_reduction_3}. In \cite{port_hamiltonian_reduction_2}, without requiring the symplecticity of the reduced basis, the gradient of the Hamiltonian vector field is approximated using a projection matrix $W$, i.e., $\nabla_y H(Uz)\approx W\nabla_{z}H_{RB}(z)$, which results in a non-canonical symplectic reduced form. The stability of the reduced model is then achieved by preserving the passivity of the original formulation. A drawback of such an approach is that, while viable for nondissipative formulations, it does not guarantee the same energy distribution of \eqref{eq:dissipative_Hamiltonian} between dissipative and null energy contributors. In the following, we show that the techniques based on symplectic geometry introduced in the previous sections can still be used in the dissipative framework described in \eqref{eq:dissipative_Hamiltonian} with limited modifications to obtain consistent and structured reduced models.  Let us consider an ortho-symplectic basis $A\in\mathbb{S}(2k,2n)$ and the reduced basis representation $y\approx Az$,  with $z=(r,s)\in\mathbb{R}^{2k}$ being the reduced coefficients of the representation and $r,s\in\mathbb{R}^{k}$ being the generalized phase coordinates of the reduced model. 
The basis $A$ can be represented as 
\begin{equation} 
A = 
\begin{bmatrix}
\enspace A_{qr} &\enspace A_{qs} \enspace \\
\enspace A_{pr} & \enspace A_{ps} \enspace
\end{bmatrix},
\end{equation}
with $A_{qr},A_{qs},A_{pr},A_{ps}\in\mathbb{R}^{n\times k}$ being the blocks, the indices of which are chosen to represent the interactions between the generalized phase coordinates of the two models, such that $q=A_{qr}r+A_{qs}s$ and $p=A_{pr}r+A_{ps}s$. Following \cite{peng_dissipative}, the symplectic Galerkin projection of \eqref{eq:dissipative_Hamiltonian} reads
\begin{equation}\label{eq:projected_dissipative}
\dot{z} = A^{+} ( \mathcal{X}_H(Az) + \mathcal{X}_F(Az) )
= \mathbb{J}_{2k}\nabla_z H_{RB}(z) + A^{+} \mathcal{X}_F(Az)
= \mathcal{X}_{H_{RB}} + A^{+} \mathcal{X}_F,
\end{equation}
with
\begin{equation}
A^{+}\mathcal{X}_F = 
\begin{bmatrix}
\enspace A_{ps}^{\top} &\enspace -A_{qs}^{\top} \enspace \\
\enspace -A_{pr}^{\top} & \enspace A_{qr}^{\top} \enspace
\end{bmatrix}
\begin{bmatrix}
\enspace 0_{n} \enspace\\
\enspace f_H \enspace
\end{bmatrix}
= 
\begin{bmatrix}
\enspace -A_{qs}^{\top} f_H \enspace \\
\enspace A_{qr}^{\top} f_H \enspace
\end{bmatrix}.
\end{equation}
We note that, in \eqref{eq:projected_dissipative}, the reduced dynamics is described as the sum of a Hamiltonian vector field and a term that, for a general choice of the symplectic basis $A$ and hence of $A_{qs}^{\top}$, does not represent a dissipative term in the form of a vertical velocity field.
The Cotangent Lift method, described in section \ref{sec:svd_based_basis_generators}, enforces by construction the structure of vertical velocity field because $A_{qs}=0$.
It can be shown \cite{peng_dissipative} that dissipativity is also preserved since the rate of energy variation of the reduced system is non-positive, i.e.,
\begin{equation}
\begin{aligned}
\nabla_s H_{RB}(Az) (A_{qr}^{\top} f_H) &= \dot{r}^{\top} (A_{qr}^{\top} f_H) 
%= -\dot{r}^{\top} A_{qr}^{\top} R(A_{qr}s)A_{qr}\dot{r} 
= - (A_{qr}\dot{r})^{\top} R(A_{qr}s) (A_{qr}\dot{r}) \leq 0.
\end{aligned}
\end{equation}
However, time discretization of the reduced dissipative model is not trivial. Even though the dissipative Hamiltonian structure is preserved by the reduction process, standard numerical integrators do not preserve the same structure at the fully discrete level.\\
A completely different approach is proposed in \cite{babak_dissipative}, where \eqref{eq:dissipative_Hamiltonian} is paired with a canonical heat bath, absorbing the energy leakage and expanding the system to the canonical Hamiltonian structure.
Consider a dissipative system characterized by the quadratic Hamiltonian $H(y)=\dfrac{1}{2}y^{\top}K^{\top}Ky$. Following \cite{TDD_formulation}, such a system admits a time dispersive and dissipative (TDD) formulation
\begin{equation}\label{eq:TDD_formulation}
\begin{cases}
\dot{y} = \mathbb{J}_{2n} K^{\top} f(t), \\
y(0)=y_0,
\end{cases}
\end{equation}
with $f(t)$ being the solution to the integral equation
\begin{equation}\label{eq:material_relation}
f(t)+\int_0^{t}\chi (t-s)f(s)ds = Ky, 
\end{equation}
also known as a \emph{generalized material relation}. The square time-dependent matrix $\chi\in\mathbb{R}^{2n\times 2n}$ is the \emph{generalized susceptibility} of the system, and it is bounded with respect to the Frobenius norm. Physically, it encodes the accumulation of the dissipation effect in time, starting from the initial condition. When $\chi=0_{2n}$, \eqref{eq:TDD_formulation} is equivalent to \eqref{eq:hamilton_equation_tbr}. Under physically natural assumptions on $\chi$ (see {\cite[Theorem 1.1, page 975]{TDD_formulation}} for more details), system \eqref{eq:TDD_formulation} admits a quadratic Hamiltonian extension (QHE) to a canonical Hamiltonian system. This extension is obtained by defining an isometric injection $I:\mathbb{R}^{2n}\mapsto \mathbb{R}^{2n}\times \mathcal{H}^{2n}$, where $\mathcal{H}^{2n}$ is a suitable Hilbert space, and reads
\begin{equation}\label{eq:embedded_TDD_system}
\begin{cases}
\dot{y}=\mathbb{J}_{2n}K^{\top} f(t), \\
\partial_t \phi = \theta (t,x),\\
\partial_t \theta = \partial_x^{2} \phi(t,x) + \sqrt{2}\delta_0(x)\cdot \sqrt{\chi} f(t),
\end{cases}
\end{equation}
where $\phi$ and $\theta$ are vector-valued functions in $\mathcal{H}^{2n}$, $\delta_0$ is the Dirac-delta function, and $f$ solves
\begin{equation*}
f(t)+\sqrt{2}\cdot\sqrt{\chi} \phi(t,0) = Ky(t).
\end{equation*}
It can be shown that system \eqref{eq:embedded_TDD_system} has the form of a conserved Hamiltonian system with the extended Hamiltonian
\begin{equation*}
H_{ex}(y,\phi,\theta) = \dfrac{1}{2}\Big ( \| Ky-\phi(t,0)\|_2^{2} +  \| \theta(t) \|^{2}_{\mathcal{H}^{2n}} + \| \partial_x \phi(t) \|^2_{\mathcal{H}^{2n}}  \Big ),
\end{equation*}
and can be reduced, while preserving its geometric structure, using any of the standard symplectic techniques.
We refer the reader to \cite{babak_dissipative} for a formal derivation of the reduced model obtained by projecting \eqref{eq:embedded_TDD_system} on a symplectic subspace and for its efficient time integration.
The method extends trivially to more general Hamiltonian functions, as long as the dissipation is linear in \eqref{eq:material_relation}.
\subsection{Non-canonical Hamiltonian systems}
The canonical Hamiltonian problem \eqref{eq:hamilton_equation_tbr} has been defined under the assumption that a canonical system of coordinates for the symplectic solution manifold is given, and the Hamiltonian vector can be represented as \eqref{eq:local_vector_field}. However, many Hamiltonian systems, such as the KdV and Burgers equations, are naturally formulated in terms of a non-canonical basis, resulting in the following description of their dynamics:
\begin{equation}\label{eq:noncanonical_Hamiltonian}
\begin{cases}
\dot{y}(t) = J_{2n}\nabla_y H(y(t)),\\
y(0) = y_0,
\end{cases}
\end{equation} 
with $J_{2n}\in\mathbb{R}^{2n\times 2n}$ being invertible and skew-symmetric. A reduction strategy, involving the non-canonical formulation \eqref{eq:noncanonical_Hamiltonian} and based on POD, has been proposed in \cite{wang_reduction}. Consider the RB ansatz $y\approx Uz$, with $U\in\mathbb{R}^{2n\times k}$ as an orthonormal basis obtained by applying the POD algorithm to the matrix of snapshots collected by solving the full model. The Galerkin projection of \eqref{eq:noncanonical_Hamiltonian} reads 
\begin{equation}\label{eq:projection_noncanonical}
\dot{z} = U^{\top} J_{2n} \nabla_{y}H(Uz),
\end{equation}
with the time derivate of the Hamiltonian function, evaluated at the reduced state, given by
\begin{equation}\label{eq:decay_projection_noncanonical}
\dot{H}(Uz) = \dot{z}^{\top} (\nabla_{z}H(Uz)) = ( \nabla_yH(Uz))^{\top} J_{2n}^{\top} UU^{\top} \nabla_yH(Uz).
\end{equation}
As expected, the Hamiltonian structure is lost in \eqref{eq:projection_noncanonical} and the energy of the system, represented by the Hamiltonian, is no longer preserved in time because $J_{2n}UU^{\top}$ is not skew-symmetric. Both issues are solved in \cite{wang_reduction} by considering a matrix $W$, with the same properties of $\mathbb{J}_{2n}$, such that the relation
\begin{equation}\label{eq:condition_preservation_noncanonical}
U^{\top}J_{2n} = W U^{\top}
\end{equation}
is satisfied. We stress that a condition similar to \eqref{eq:condition_preservation_noncanonical}  naturally holds in the canonical Hamiltonian setting for a symplectic basis and has been used to derive Hamiltonian reduced models using the symplectic Galerkin projection. A candidate $W$ is identified in \cite{wang_reduction} by solving the normal equation related to \eqref{eq:condition_preservation_noncanonical}, i.e.
$W = U^{\top} J_{2n} U$.
For invertible skew-symmetric operators $J_{2n}$ that might depend on the state variables $y$, Miyatake has introduced in \cite{noncanonical_DEIM} an
hyper-reduction technique that preserves the skew-symmetric structure of the $J_{2n}$ operator. \\
Formulation \eqref{eq:noncanonical_Hamiltonian} is further generalized with the characterization of the phase-space as a Poisson manifold, defined as a $2n_P$-dimensional differentiable manifold $\mathcal{M}_{P}$ equipped with a Poisson bracket $\{ \cdot , \cdot \}:C^{\infty}(\mathcal{M}_{p})\times C^{\infty}(\mathcal{M}_{p}) \mapsto C^{\infty}(\mathcal{M}_{p}) $ satisfying the conditions of bilinearity, skew-symmetry, the Jacobi identity, and the Leibniz' rule. Since derivations on $C^{\infty}(\mathcal{M}_{P})$ are represented by smooth vector fields, for each Hamiltonian function $H\in C^{\infty}(\mathcal{M}_{P})$,
there exists a vector $\mathcal{X}_H$ that determines the following dynamics,
%there exists a vector $\mathcal{X}_H(g):=\{ H,g \}$, known as the Hamiltonian vector field, determining the following Hamiltonian dynamics
\begin{equation}\label{eq:noncanonical_statedependent_Hamiltonian}
\begin{cases}
\dot{y}(t) = \mathcal{X}_H(y) = J_{2n_P}(y)\nabla_y H(y(t)),\\
y(0) = y_0,
\end{cases}
\end{equation}
with the Poisson tensor $J_{2n_p}$ being skew-symmetric, state-dependent, and generally not invertible.
The flow of the Hamiltonian vector field $\mathcal{X}_H(y)$, which is a Poisson map and therefore preserves the Poisson bracket structure via its pullback, also preserves the rank $2n$ of the Poisson tensor $J_{2n_P}(y)$.
 Moreover, $r=2n_{P}-2n$ represents the number of independent nonconstant functions on $\mathcal{M}_P$ that $\{ \cdot, \cdot \}$ commutes with all the other functions in $C^{\infty}(\mathcal{M}_P)$. These functions are known as Casimirs of the Poisson bracket and their gradients belong to the kernel of $J_{2n_P}(y)$, making them independent of the dynamics of \eqref{eq:noncanonical_statedependent_Hamiltonian} and only representing geometric constraints on configurations of the generalized phase-state space. \\
An interesting relation between symplectic and Poisson manifolds is offered by the Lie-Weinstein splitting theorem, stating that locally, in the neighborhood $\mathcal{U}_{y^{*}}$ of any point $y^{*}\in\mathcal{M}_P$, 
a Poisson manifold can be split into a $2n$-dimensional symplectic manifold $\mathcal{M}$ and an $r$-dimensional Poisson manifold $M$. Following on this result, Darboux' theorem guarantees the existence of local coordinates $(q_1,\dots,q_n,p_1,\dots,p_n,c_1,\dots, c_r)$, where $\{ q_i,p_i \}_{i=1}^{2}$ corresponds to canonical symplectic coordinates and $\{ c_i \}_{i=1}^{r}$ are the Casimirs, such that the Poisson tensor $J_{2n_p}(y)$ is recast, via Darboux' map, in the canonical form $J_{2n_p}^{C}$, i.e.,
\begin{equation*}
J_{2n_p}^{C} =
    \begin{blockarray}{ccc}
     \enspace 2n & r & \\
      \begin{block}{[cc]c}
      \mathbb{J}_{2n} & 0 & 2n \\ 
      0 & 0 & r \\
      \end{block}
    \end{blockarray},
\end{equation*}
with $\mathbb{J}_{2n}\in\mathbb{R}^{2n\times 2n}$ being the canonical Poisson tensor defined in \ref{eq:Poisson_tensor}.\\
In \cite{pagliantini_2}, a quasi-structure-preserving algorithm for problems of the form \eqref{eq:noncanonical_statedependent_Hamiltonian} has been proposed, leveraging the Lie-Weinstein splitting, an approximation of the Darboux' map and traditional symplectic RB techniques.  
Let
\begin{equation}\label{eq:noncanonical_statedependent_Hamiltonian_discrete}
\begin{cases}
y^{j+1} = y^{j} + \Delta t J_{2n_p}(\tilde{y}^{j})\nabla_y H(\tilde{y}^{j}), \\
y^{0} = y_0,
\end{cases}
\end{equation}
be the fully-discrete formulation of \eqref{eq:noncanonical_statedependent_Hamiltonian}, where $j$ is the integration index, and $\tilde{y}^{j}$ represents intermediate state/states dictated by the temporal integrator of choice.
Given $\mathcal{M}_{P,j}$, an open subset of $\mathcal{M}_P$ comprising the discrete states $y^{j}$, $\tilde{y}^{j}$, and $y^{j+1}$, the authors of \cite{pagliantini_2} introduce an approximation $\varphi_{j+\frac{1}{2}}:\mathcal{M}_{P,j}\mapsto \mathcal{M}_s \times \mathcal{N}_j$ of the Darboux' map at $\tilde{y}^{j}$, with $\mathcal{M}_s$ being a $2N$-dimensional canonical symplectic manifold and $\mathcal{N}_j$ approximating the null space of the Poisson structure. The proposed approximation exploits a Cholesky-like decomposition (see \cite[Proposition 2.11, page 1708]{pagliantini_2}) of the noncanonical rank-deficient $J_{2n_p}(\tilde{y}^{j})$ and exactly preserves the dimension of $\mathcal{N}_j$, hence the number of independent Casimirs. By introducing the natural transition map $T_j:=\varphi_{j+\frac{1}{2}}\cdot \varphi_{j-\frac{1}{2}}^{-1}$ between the neighboring and overlapping subsets $\mathcal{M}_{j-1}$ and  $\mathcal{M}_{j}$, problem \eqref{eq:noncanonical_statedependent_Hamiltonian_discrete} is locally recast in the canonical form 
\begin{equation}\label{eq:noncanonical_statedependent_Hamiltonian_discrete_local}
\begin{cases}
\bar{y}^{j+1} = T_j \bar{y}^{j} + \Delta t J_{2n_p}^{C}\nabla_{\bar{y}} H^j(\bar{\tilde{y}}^{j}), \\
\bar{y}^{0} = y_0,
\end{cases}
\end{equation}
where $\bar{y}^{j+1}:=\varphi_{j+\frac{1}{2}}y^{j+1}$, $\bar{y}^{j}:=\varphi_{j+\frac{1}{2}}y^{j}$, $\bar{\tilde{y}}^{j+1}:=\varphi_{j+\frac{1}{2}}\tilde{y}^{j}$, and $H^{j}(\bar{\tilde{y}}^{j}):=H(\varphi_{j+\frac{1}{2}}^{-1}(\bar{\tilde{y}}^{j}))$.
Even though the flow of \eqref{eq:noncanonical_statedependent_Hamiltonian_discrete_local} is not a \emph{global} $J_{2n_p}^{C}$-Poisson map because the splitting is not exact, the approximation is \emph{locally} structure-preserving for each neighborhood $\mathcal{M}_{P,j}$. By exploiting a similar splitting principle, the canonical Poisson manifold $\mathcal{M}_s\times \mathcal{N}_j$ is projected on a reduced Poisson manifold $\mathcal{A}\times \mathcal{N}_j$, with the reduction acting only on the symplectic component of the splitting and $\text{dim}(\mathcal{A})=2k\ll 2n$. The corresponding reduced model is obtained via Galerkin projection of \eqref{eq:noncanonical_statedependent_Hamiltonian_discrete_local} using an orthogonal $J_{2k}^{C}$-symplectic basis of dimension $2k$, generated via a greedy iterative process inspired by the symplectic greedy method described in Section \ref{sec:greedy_approach}. Different theoretical estimates and numerical investigations show the proposed technique's accuracy, robustness, and conservation properties, up to errors in the Poisson tensor approximation. 

\section{Conclusion}\label{sec:conclusion}
We provided an overview of model reduction methods for Hamiltonian problems.
The symplectic Galerkin projection has been discussed as a tool to generate a reduced Hamiltonian approximation of the original dynamics. PSD algorithms used to compute low-order projection on symplectic spaces have been introduced and compared. Such strategies have been classified in ortho-symplectic and symplectic procedures, depending on the structure of the computed RB. A greedy alternative for the generation of ortho-symplectic basis, characterized by an exponentially fast convergence, has been illustrated as an efficient iterative approach to overcome the computational cost associated with SVD-based techniques that require a fine sampling of the solution manifold of the high-dimensional problem. 
The potential local low-rank nature of Hamiltonian dynamics has been addressed by a symplectic dynamical RB method. The innovative idea of the dynamical approach consists in evolving the approximating symplectic reduced space in time along a trajectory locally constrained on the tangent space of the high-dimensional dynamics.
For problems where the Hamiltonian dynamics is coupled with a dissipative term, structure-preserving reduced models can be constructed with the symplectic reduction process by resorting to an extended non-dissipative Hamiltonian reformulation of the system. Finally, we have described  RB strategies to reduce problems having a non-canonical Hamiltonian structure that either enforce properties typical of a symplectic basis or use canonical symplectic reductions as an intermediate step to preserve the structure of the original model.

\bibliographystyle{unsrt}  
%\bibliography{references}  %%% Remove comment to use the external .bib file (using bibtex).
%%% and comment out the ``thebibliography'' section.

%%% Comment out this section when you \bibliography{references} is enabled.

\end{document}